%% file: 3d_whitney220716.tex
\documentclass[11pt]{article}

\input{minimalJ}

\usepackage{verbatim}
\usepackage{enumitem}
\usepackage[all]{xy}
\usepackage{subfig}

\setenumerate{label={\textnormal{(\arabic*)}}}

\makeatletter
\newcommand{\dotminus}{\mathbin{\text{\@dotminus}}}

\newcommand{\@dotminus}{%
  \ooalign{\hidewidth\raise1.2ex\hbox{.}\hidewidth\cr$\m@th-$\cr}%
}
\makeatother

\makeatletter
\newcommand{\oset}[3][0ex]{%
  \mathrel{\mathop{#3}\limits^{
    \vbox to#1{\kern-0.1\ex@
    \hbox{$\scriptstyle#2$}\vss}}}}
\makeatother

\theoremstyle{plain}

\usepackage{xcolor}
\colorlet{darkishRed}{red!80!black}
\colorlet{darkishBlue}{blue!60!black}
\colorlet{darkishGreen}{green!60!black}
\hypersetup{
    draft = false,
    bookmarksopen=true,
    colorlinks,
    linkcolor={red!60!black},
    citecolor={green!60!black},
    urlcolor={blue!60!black}
}

\title{Dual matroids of 2-complexes -- revisited}

\author{Johannes Carmesin\footnote{Funded by EPSRC, grant number EP/T016221/1}\\
  {University of Birmingham}
}

\DeclareMathOperator{\Sbb}{\mathbb{S}}

\newcommand{\Sthree}{$\Sbb^3$}

\mcm{\Fbb}{0}{\mathbb{F}}

\usepackage{xr}
\begin{document}
\maketitle

\begin{abstract}
We prove that simply connected local 2-dimensional simplicial complexes embed in 3-space if and only if their dual matroids are graphic.
Examples are provided that the assumptions of simply connectedness and locality are necessary.
This may be regarded as a 3-dimensional analogue of Whitney's planarity criterion from 1932.
\end{abstract}

\section{Introduction}

A fundamental theorem at the interface of graph theory and matroid theory is Whitney's planarity criterion: a graph is planar if and only if its dual matroid is graphic \cite{Whitney32}. In this paper we prove a 3-dimensional analogue of this theorem.

Since every graph has an embedding in 3-space, it seems natural to consider embeddings of 2-dimensional generalisations of graphs; that is, embeddings of 2-dimensional simplicial complexes in 3-dimensional space. Similar to the situation in the plane, every 2-dimensional simplicial complex $C$ embedded in 3-space has a dual graph: its vertices are the chambers of the embedding and two of these vertices are joined by an edge if they share a face of $C$. As for the plane, this dual graph can be computed from a matroid, the matroid represented by the edge/face incidence matrix of $C$. Wouldn't it be wonderful if a 2-dimensional simplicial complex embedded in 3-space if and only if its dual matroid was graphic? It turns out that we have to be a bit more careful for the following reasons.

\begin{figure} [htpb]
\begin{center}
   	  \includegraphics[height=2cm]{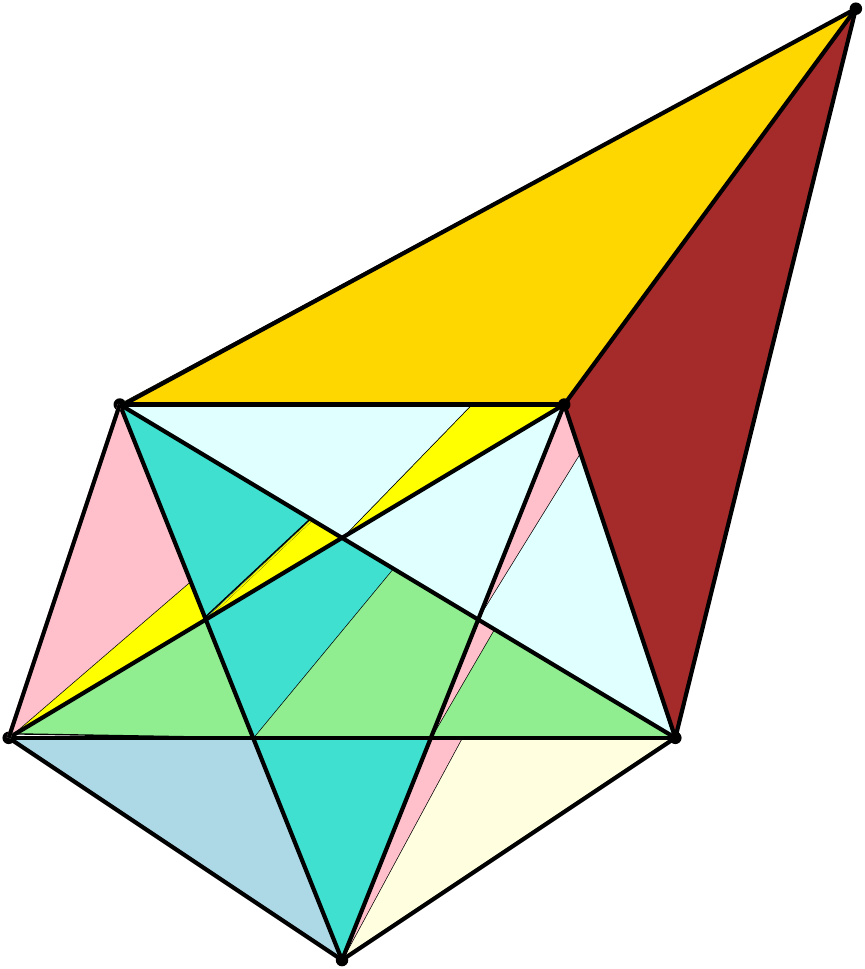}
   	  \caption{The cone over the graph $K_5$ has a dual matroid that is a sum of loops, yet it does not embed in any 3-manifold analogous to the fact that $K_5$ does not embed in the plane.
   	  }\label{fig:cone}
\end{center}
   \end{figure}

\begin{enumerate}
 \item Triangulations of Poincar\'e-spheres have graphic dual matroids but no embeddings in $\Sbb^3$ in general, as matroids work on the level of homology but are too rough to detect differences on the level of homotopy, roughly speaking. Hence it seems natural to restrict our attention to simply connected simplicial complexes $C$;
 \item The cone over the graph $K_5$ does not embed in any 3-manifold yet has a trivial graphic dual matroid. The notion of \lq locality\rq\ has been tailored to prevent this type of obstruction \cite{3space4}, see \autoref{fig:cone} and \autoref{sec2} for a definition;
 \item While for embeddable $C$, its edge/face incidence matrix represents a regular matroid and so this matroid can be defined independently of the field $k$ we work over, for $C$  generic this is not the case.
\end{enumerate}

Under these three necessary conditions we answer the above question affirmatively.

\begin{thm}\label{intro-main}
For every field $k$, a $k$-local simply connected 2-dimensional simplicial complex $C$ embeds in 3-space if and only if its $k$-dual matroid is graphic; in this case $k$-dual matroids over different fields $k$ are isomorphic.
\end{thm}

A 3-dimensional analogue of MacLane's planarity criterion follows from \autoref{intro-main}, see \autoref{concluding}.

{\bf Related results.}
The special case of \autoref{intro-main} for fields $\Fbb_p$ with $p\neq 2$ was proved by the author 5 years ago \cite{3space4}.
That paper is part of a series of five papers of the author in which they solve questions asked by Lov\'asz, Pardon and Wagner by proving a 3-dimensional analogue of Kuratowski's theorem, which yielded a polynomial algorithm to check embeddability for simply connected 2-dimensional simplicial complexes. Follow-up works include an analogue of Whitney's Uniqueness theorem by Georgakopoulos and Kim \cite{georgakopoulos20212}, an explicit characterisation of outerspacial 2-complexes \cite{outerspace}, and a Heawood-type theorem for 2-complexes \cite{3dheawood}.

The remainder of the paper is structured as follows. In \autoref{sec2} we explain the basic tools from 3-dimensional combinatorics required in this paper.
The next two sections are concerned with constructing embeddings in 3-space for local 2-dimensional simplicial complexes $C$ with graphic dual matroids, one implication of \autoref{intro-main}.
In \autoref{sec3} we construct a rotation framework, some combinatorial data that will help embed $C$, and in \autoref{sec4} we show that in the presence of a graphic dual matroid, locality and simply connectedness, this rotation framework has some stronger properties, which will then allow us to construct an embedding of $C$ in a 3-manifold; by Perelman's theorem this 3-manifold must be $\Sbb^3$. In \autoref{concluding} we deduce a 3D MacLane theorem form \autoref{intro-main}.

For matroid background we refer the reader to \cite{Oxley2}, for topology use \cite{Hatcher} and for embeddings of graphs use \cite{MoharThomassen}.

\section{Basics from 3-dimensional Combinatorics}\label{sec2}

In this paper we abbreviate 2-dimensional simplicial complex simply by \emph{simplicial complex}.
We denote simplicial complexes by triples  $(V,E,F)$, where $V$ is its set of vertices, $E$ is its set of edges and $F$ is its set of faces. We follow the convention that every edge of a simplicial complex is incident with a face. Given a simplicial complex $C$ and a vertex $v$ of $C$, the \emph{link graph} $L(v)$ is the graph whose vertices are the edges incident with $v$ and two of them are adjacent if they are in a common face; that is, its edge set is the set of faces containing $v$. In this paper we suppress from our notation an injection from the vertices and edges of the link graph to $E$ and $F$, respectively; and simply consider the vertices of the link graph as edges of $C$, for instance.

A \emph{directed} simplicial complex is a simplicial complex together with a choice of direction at each edge and a choice of orientation at each face.
The \emph{incidence vector} of an edge $e$ is the vector in the face-space $\{0,\pm 1\}^F$, where the vector has the entry $0$ in coordindate $f$ if $e$ and $f$ are not incident, $1$ if the direction at $e$ is positive in the orientation of $f$ and $-1$ otherwise.
Throughout the paper we fix a field $k$.
Given $k$ and a directed simplicial complex $C$, its \emph{$k$-dual matroid} is the $k$-vector matroid whose ground set is the
set $F$ of faces of $C$ and its circuit space is generated by the incidence vectors of the edges.
Occasionally, we shall omit the explicit dependence on the field $k$ and assume that it is implicitly given, for example we may abbreviate \lq $k$-dual matroid\rq\ by \lq dual matroid\rq. We denote the dual matroid of $C$ by $M^*(C)$.
\begin{rem}
For different directed simplicial complexes with the same underlying simplicial complex, dual matroids are isomorphic. Hence we shall sometimes omit the term \lq directed\rq\ and assume that some direction is implicitly given.
\end{rem}

One implication of \autoref{intro-main} is given by the following.

\begin{thm}[{\cite[\autoref{2combi_intro}]{3space2}}]\label{x1}
Let $C$ be a (directed) 2-dimensional simplicial complex embedded into \Sthree. Then the edge/face
incidence matrix of $C$
represents over the integers (and hence over any field $k$) a
matroid $M=M^*(C)$ which is equal to the cycle matroid of the dual graph of the embedding.
\end{thm}

We say that a simplicial complex $C$ is \emph{$k$-local} if for every vertex $v$ of $C$ the dual of the cycle matroid of the link graph $L(v)$ is equal to the restriction of the $k$-dual matroid of $C$ to those faces that contain $v$; in formulas: $M(L(v))^*=M^*(C)\restric_E$, where $E=E(M(L(v)))$.

The goal of this paper is to go the other way than \autoref{x1}: to prove that every $k$-local simply connected simplicial complex whose $k$-dual matroid is graphic embeds in 3-space.
Towards this, in \cite{3space4}, we proved\footnote{The argument is a combination of Sublemmas 3.6 and 3.7 in that paper and given in the environment marked as \lq Proof of Theorem 1.2\r.} the following. A simplicial complex is \emph{locally 2-connected} if all its link graphs are 2-connected.

\begin{lem}\label{red_to_2con}
 Assume that every simply connected locally 2-connected $k$-local simplicial complex whose dual matroid is graphic embeds in 3-space.

 Then every simply connected $k$-local simplicial complex whose dual matroid is graphic embeds in 3-space.
\end{lem}

Our aim is to prove the following.

\begin{lem}\label{rest}
Every simply connected locally 2-connected $k$-local simplicial complex whose dual matroid is graphic embeds in 3-space.
\end{lem}

\begin{proof}[Proof that \autoref{rest} implies \autoref{intro-main}.]
By \autoref{x1} every simplicial complex embedded in 3-space has a graphic dual matroid; and these dual matroids are isomorphic over different fields.
Conversely, by \autoref{rest} and \autoref{red_to_2con}  every simply connected $k$-local simplicial complex whose $k$-dual matroid is graphic  embeds in 3-space.
\end{proof}

The proof of \autoref{rest} takes the rest of the paper. First we do some preparation.
A \emph{generalised simplicial complex} is a simplicial complex plus some parallel faces; that is, we add copies of faces with the same incidences.

\begin{rem}
The definitions for simplicial complexes extend to generalised simplicial complexes in the natural way, and we will not always point this out explicitly.
\end{rem}

We shall make use of the following standard facts from planar graph theory.
A \emph{simultaneous dual embedding} of two graphs $G$ and $H$ with a bijection between their edge sets is a pair of embeddings, one for $G$ and one for $H$ such that these embeddings are disjoint except that for every edge $e$ its embeddings for $G$ and $H$ intersect in the midpoint of that edge $e$.

\begin{thm}[Whitney's planarity criterion (slight variation\footnote{For 3-connected graphs this statement follows directly from Whitney's planarity criterion in the form of \cite{DiestelBookCurrent}. The generalisation to 2-connected graphs easily follows from the 2-separator theorem.})]\label{min_whitney}
 Let $G$ be a 2-connected graph such that there is a graph $H$ whose cycle matroid is dual to the cycle matroid of $G$. Then $G$ and $H$ have a simultaneous dual plane embedding.
\end{thm}

\begin{lem}[\cite{MoharThomassen}]\label{lem_2con}
 In a 2-connected plane graph every face is bounded by a cycle.
\end{lem}

A \emph{rotator} at a vertex of a graph $G$ is a cyclic orientation of its incident edges.
An embedding of a graph in the 2-sphere is combinatorially characterised through a \emph{rotation system}; that is, a choice of rotator at each of its vertices.
We say that a rotation system is \emph{induced} by an embedding in $\Sbb^2$ (with a fixed orientation) if for every vertex the cyclic orientation of the incident edges when we go around the vertex in the chosen orientation of $\Sbb^2$ is equal to its rotator. Rotation systems are a useful tool to study embeddings of graphs in orientable surfaces and there is an unoriented analogue \cite{MoharThomassen}. There are similar notions for simplicial complexes and embeddings in 3-space, as follows (see \cite{3space2} for details on alternative notions).

Given a simplicial complex $C$, a \emph{rotation framework} of $C$ consists of choices of embeddings of the link graphs at all vertices of $C$ such that these embeddings are compatible in that at every edge $e=vw$ of $C$, the chosen rotators in $L(v)$ and $L(w)$ at $v$ are reverse or agree.
We colour an edge \emph{green} if the rotators are reverse, otherwise we colour it red.
A rotation framework is \emph{even} if the number of red edges is even at every cycle of (the 1-skeleton of) $C$.
We say that a rotation framework $\Sigma$ is \emph{induced} by an embedding of a simplicial complex $C$ in a 3-manifold if at each vertex $v$ of $C$, there is a (small) 2-sphere around $v$ in the 3-manifold such that $C$ intersects that sphere in the link graph at $v$ and the embedding of the link graph on the outside of that 2-sphere is equal to its embedding for $\Sigma$.

\begin{thm}[Skopenkov 1994 \cite{Skopenkov94}]\label{skop}
 A simplicial complex that admits an even rotation framework $\Sigma$ embeds in some compact 3-manifold such that the rotation framework induced by the embedding is equal to $\Sigma$.
\end{thm}

\begin{eg}
A torus obstruction, see \cite{3space5}, has a rotation framework but not an even one.
\end{eg}

\section{Constructing rotation frameworks}\label{sec3}

In this section we prove \autoref{construct_PRS}, which is used in the proof of \autoref{intro-main}.

Let $C$ be a simplicial complex such that its $k$-dual matroid $M$ is graphic and let $G$ be a graph whose cycle matroid is $M$. We refer to these conditions as \emph{a simplicial complex $C$ with dual graph $G$}. Given a vertex $v$ of $C$, denote by $G_v$ the graph obtained from $G$ by deleting all edges whose corresponding faces of $C$ at not incident with $v$.

\begin{lem}\label{construct_PRS}
Let $C$ be a locally 2-connected $k$-local generalised simplicial complex with dual graph $G$. Then $C$ has a rotation framework such that for every vertex $v$ the dual graph of the embedding of the link graph $L(v)$ at $v$ is equal to the graph $G_v$.
\end{lem}

\begin{proof}
First we construct a rotation framework.
For a vertex $v$ of $C$, denote by $M_v$ the dual matroid of the cycle matroid of the link graph $L(v)$.

\begin{sublem}\label{dual}
 The link graph $L(v)$ has a plane embedding whose dual graph is $G_v$.
\end{sublem}

\begin{proof}
By locality, the cycle matroid of the graph $G_v$ is equal to the matroid $M_v$. Hence by Whitney's planarity criterion, the matroid $M_v$ is the cycle matroid of a planar graph; that is, the graph $G_v$ is planar. By 2-connectedness of $L(v)$, the graphs $L(v)$ and $G_v$ have a simultaneous dual plane embedding, see \autoref{min_whitney}.
\end{proof}

Let $e$ be an edge of $C$.

\begin{sublem}\label{choice}
The faces incident with $e$ in $C$ form a cycle in the graph $G$.
\end{sublem}
\begin{proof}
Let $v$ be an endvertex of the edge $e$. Through \autoref{dual} pick an embedding of $L(v)$ with dual graph $G_v$.
The edges incident with $e$ in $L(v)$ form a face boundary in the embedding of $G_v$, and by 2-connectedness, see \autoref{lem_2con}, this face boundary is a cycle.
The edges incident with $e$ in $L(v)$ are precisely the faces of $C$ that are incident with $e$. So this set of faces forms a cycle in the graph $G$.
\end{proof}

By \autoref{dual} we can pick an embedding of each link graph $L(v)$ such that its dual graph is $G_v$. By \autoref{choice} the rotator at some vertex $e$ of $L(v)$ is determined by the cyclic ordering of the edges of a cycle. Hence when we consider $e$ as a vertex of the other link graph containing it, in there the rotator is the same or its reverse. Thus these choices define a rotation framework. This completes the proof.
\end{proof}

We say that a rotation framework as in \autoref{construct_PRS} is \emph{induced} by $G$.

\section{Controlling rotation frameworks through locality}\label{sec4}

In this section we prove \autoref{rest}, which implies \autoref{intro-main} as shown above.
We start with some preparation.

Let $C$ be a simplicial complex with dual graph $G$. Given a vertex $v$ of $G$ and an edge $e$ of $C$, we say that a face of $C$ is \emph{$(v,e)$-incident} if it is incident with the edge $e$ of $C$ and when considered as an edge of $G$ it is incident with the vertex $v$.
We say that the pair $(C,G)$ is \emph{sparse} if for every pair $(v,e)$ with $v\in V(G)$ and $e\in E(C)$, there are exactly two or exactly zero faces of $C$ that are
$(v,e)$-incident.

\begin{lem}\label{max2}
Let $C$ be a locally 2-connected $k$-local generalised simplicial complex $C$ with dual graph $G$. Then $(C,G)$ is sparse.
\end{lem}

\begin{eg}
 If $C$ as in \autoref{max2} is embedded in $\Sbb^3$ with dual graph $G$ (for example $C$ is a triangulation of the 2-sphere and $G$ is a parallel graph), \autoref{max2} implies that the boundary of every chamber is a surface.
\end{eg}

\begin{proof}[Proof of \autoref{max2}]
Take a rotation framework $\Sigma$ of $C$ induced by $G$ via \autoref{construct_PRS}.
Let a vertex $v$ of $G$ and an edge $e$ of $C$ be given. Assume that there is some face $f$ of $C$ that is incident with $e$ in $C$ and incident with $v$ when considered as an edge of $G$.

In $C$ pick an endvertex of the edge $e$ arbitrarily, call it $x$. Consider the embedding of the link graph $L(x)$ at $x$ given by the rotation framework $\Sigma$.
In there, $f$ is an edge.
By locality, the faces of the embedding of the link graph $L(x)$ that have the edge $f$ in their boundary correspond to the vertices of $G$ that are incident with $f$ when considered as an edge of $G$. Denote the face of the embedding of $L(x)$ corresponding to the vertex $v$ by $v'$. By assumption the link graph $L(x)$ is 2-connected, so the face $v'$ is bounded by a cycle by \autoref{lem_2con}.
So in the link graph $L(x)$ there is exactly one edge that is incident with the vertex $e$ and the face $v'$ aside from $f$.
So by locality there is exactly one edge of $G$ aside from $f$ that is incident with $v$ and incident with $e$ when considered a face of $C$. Hence $(C,G)$ is sparse.
\end{proof}

The \emph{face-degree} of an edge is the number of its incident faces.

\begin{lem}\label{kill-pre}
 Let $C$ be a locally 2-connected $k$-local generalised simplicial complex with dual graph $G$. Let $f$ be a face of $C$ such that all its incident edges have face-degree at least three.
 Let $\Sigma$ be a rotation framework of $C$ induced by $G$. Then the number of red edges (with respect to $\Sigma$) in the face boundary of $f$ is even.
\end{lem}

\begin{proof}
Let $e_1$, $e_2$ and $e_3$ be the edges in the boundary of $f$. Denote the vertices in the face $f$ by $v_1$, $v_2$ and $v_3$ where $e_i$ is incident with $v_i$ and $v_{i+1}$. Note that flipping the orientation of a plane embedding at a link graph for $\Sigma$ gives us a rotation framework that is induced by the same dual graph $G$. As all edges of $C$ incident with the face $f$ of $C$ have face-degree at least three, such a flip would change the colours at all incident edges. Thus such a flipping preserves the parity of the number of red edges of $C$ on the face $f$ of $C$.
By flipping the orientation in the embedding of the link graph at $v_2$ in $\Sigma$ if necessary, assume that $e_1$ is green. By flipping at $v_3$ if necessary assume additionally that $e_2$ is green. With this modification in place we shall in fact prove that the final edge $e_3$ of $f$ is green as well.

In the embedding of the link graph $L(v_{i+1})$ at $v_{i+1}$ (with respect to $\Sigma$) denote the edge just before the edge $f$ at the rotator at $e_i$ by $x(i,1)$ and the edge just after $f$ by $x(i,2)$. Similarly, in the embedding of $L(v_{i+1})$ denote the edge just before $f$ at the rotator at $e_{i+1}$ by $y(i,1)$ and the edge just after $f$ by $y(i,2)$, see \autoref{fig:edge_f}.

   \begin{figure} [htpb]
\begin{center}
   	  \includegraphics[height=3cm]{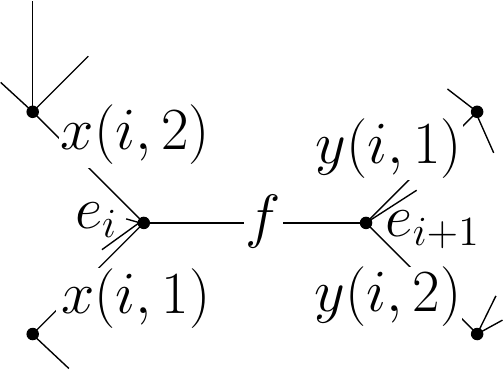}
   	  \caption{A cut-out of the link graph $L(v_{i+1})$ at the vertex $v_{i+1}$. The edge $f$ is incident with the vertices $e_i$ and $e_{i+1}$ and one face of the embedding contains the edges $x(i,1)$ and $y(i,2)$, while the other contains $x(i,2)$ and $y(i,1)$, compare \autoref{sym}. Note that $f$ is not a bridge by 2-connectedness of the link graph and so its two incident faces are distinct.}\label{fig:edge_f}
\end{center}
   \end{figure}

Considering $f$ as an edge of $G$, pick one of its endvertices arbitrarily and denote it by $b$.

\begin{sublem}\label{sym}
For $i\in \Zbb_3$, if in $G$ the vertex $b$ is incident with the edge $x(i,1)$, it is incident with  $y(i,2)$.
\end{sublem}

\begin{proof}
By \autoref{max2} in $G$ the vertex $b$ is incident with at most two edges that in $C$ correspond to faces incident with the edge $e_i$. As $e_i$ is incident with $f$ and $x(i,1)$, it cannot be incident $x(i,2)$ (which is different from $x(i,1)$ and $f$ as $e_i$ has face-degree at least three by assumption).

Note that in the embedding of $L(v_{i+1})$, one face containing $f$ includes the edges $x(i,1)$ and $y(i,2)$, while the second face containing $f$ includes the edges $x(i,2)$ and $y(i,1)$.
Via locality, consider $b$ as a unique face of the embedding of $L(v_{i+1})$. So the face $b$ is the face containing $f$ but not $x(i,2)$. This face contains the edge $y(i,2)$.
\end{proof}

By symmetry assume that in $G$ the vertex $b$ is incident with the edge $x(1,1)$.
Suppose for a contradiction that the edge $e_3$ is red. Then $y(3,2)=x(1,2)$.
Now we make the following chain of arguments:

\begin{itemize}
 \item $b$ is incident with $y(1,2)$ (by \autoref{sym});
 \item $y(1,2)=x(2,1)$ (as $e_1$ is green);
 \item $b$ is incident with $y(2,2)$ (by \autoref{sym});
  \item $y(2,2)=x(3,1)$ (as $e_1$ is green);
 \item $b$ is incident with $y(3,2)$ (by \autoref{sym});
 \item $y(3,2)=x(1,2)$ (as $e_1$ is red);
\end{itemize}

In particular, in $G$ the vertex $b$ is incident with the edges $f$, $x(1,1)$ and $x(1,2)$. As all of these faces of $C$ are incident with the edge $e_1$ of $C$, we conclude that $(C,G)$ cannot be sparse; this is a contradiction to \autoref{max2}.
So the edge $e_3$ is green. So the number of red edges in the boundary of $f$ is even.

\end{proof}

Given a generalised simplicial complex $C$ obtained from a generalised simplicial complex $C'$ by deleting faces, we say that a rotation framework $\Sigma'$ of $C'$ \emph{induces} a rotation framework $\Sigma$ of $C$ if for each vertex $v$ the embedding of the link graph at $v$ in $C$ for $\Sigma$ is obtained from the embedding of the link graph at $v$ in $C'$ for $\Sigma'$ by deleting edges.

\begin{lem}\label{junkify}
 Let $C$ be a locally 2-connected simplicial complex with a rotation framework $\Sigma$ induced by a graph $G$.
 There is a generalised simplicial complex $C'$ obtained from $C$ by adding parallel faces that admits a rotation framework $\Sigma'$ that induces $\Sigma$ and such that every edge has face-degree at least three.

 Moreover $\Sigma'$ is induced by a graph $G'$ and $G'$ is a subdivision of $G$.
\end{lem}

\begin{proof}
 Note that by local 2-connectivity, all edges have face-degree at least two. We prove the statement by induction on the number of edges of face-degree two.
 The induction starts with all edges having face-degree at least three. And we pick $C'=C$ in this case.
 So assume that $C$ has at least one edge $e$ of face-degree two. Pick a face $x$ incident with $e$. We obtain $C'$ from $C$ adding a copy $x'$ of $x$ with the same incidences.

 Now we define a rotation framework $\Sigma'$ for $C'$.
  Denote the three vertices incident with $x'$ by $v_1$, $v_2$ and $v_3$, where $e=v_3v_1$.
  We extend the embeddings of the link graphs $L(v_i)$ for $C$ by adding an edge labelled $x'$ so that $x'$ and $x$ bound a face of size two.
  We pick these embeddings inductively along the path $v_1v_2v_3$ ensuring compatibility on the two edges of that paths. Finally we observe that the edge $e$ has face-degree three in $C'$ and there is a unique cyclic ordering on a set of size three, so compatibility at $e$ is automatically satisfied. This completes the construction of the rotation framework $\Sigma'$ of $C'$.
By construction $\Sigma'$ induces $\Sigma$.

To show the \lq Moreover\rq-part, we construct a graph $G'$ from $G$ by subdividing the edge $x$ and denoting the two subdivision edges by \lq$x$\rq\ and \lq$x'$\rq. By construction the dual matroid of $C'$ is equal to the cycle matroid of $G'$. So $C'$ is a simplicial complex with dual graph $G'$.
The construction of the embeddings at the link graphs $L(v_i)$ for $C'$ for $\Sigma'$ ensures that $\Sigma'$ is induced by the graph $G'$.
\end{proof}

\begin{proof}[Proof of \autoref{rest}]
Let $C$ be a simply connected locally 2-connected $k$-local simplicial complex such that its $k$-dual matroid represents a graph $G$.
By \autoref{construct_PRS}, $C$ has a rotation framework that is induced by $G$.
We construct $C'$ from $C$ via \autoref{junkify}. By that lemma $C'$ has a rotation framework $\Sigma'$ that is induced by some graph and all its edges have face-degree at least three. Note that $C'$ is simply connected, locally 2-connected and $k$-local.
By \autoref{kill-pre}, the number of red edges (with respect to $\Sigma'$) at every face boundary of $C'$ is even.
By simply connectedness, the face boundaries generate all cycles over the field $\Fbb_2$. Hence the number of red edges on every cycle is even; that is, the rotation framework $\Sigma'$ is even. By \autoref{skop}, $C'$ has an embedding in a compact 3-manifold $X$.
As $C'$ is simply connected, also $X$ is simply connected by Van-Kampen's Theorem. So the  3-manifold $X$ is isomorphic to $\Sbb^3$ by Perelman's theorem.
By restriction, this embedding induces an embedding of $C$ in $\Sbb^3$ whose induced planar rotation system is $\Sigma$.
\end{proof}

\begin{proof}[Proof of \autoref{intro-main}]
 We have just completed the proof of \autoref{rest}. Directly below its statement, we proved that it implies \autoref{intro-main}.
\end{proof}

\section{Concluding Remark}\label{concluding}

We conclude with the following application.
Given a field $k$, the \emph{cycle space} $Z_2^k(C)$ of a simplicial complex $C$ with $F$ its set of faces is the kernel of the boundary map $\partial_2$; that is, it consists of those vectors of $k^F$ whose boundary is zero. The \emph{cycle matroid} of $C$ over $k$ is the
matroid whose cycle space is $Z_2^k(C)$.
The vector space $Z_2^k(C)$ is the orthogonal complement of the cycle space of the $k$-dual
matroid of $C$; thus the cycle matroid of $C$ is dual to the dual matroid of $C$.
Given a field $k$ and a vector space $k^F$, a family of vectors of $k^F$ is \emph{$k$-sparse} if in every coordinate  all of them have zero entries or exactly two of them have non-zero entries and these entries are $1$ and $-1$.

Whitney's characterisation of planarity implies MacLane's
characterisation of planarity once we make the following observation.

\begin{obs}\label{MC-lem}
Given a field $k$, a matroid is graphic if and only if its dual matroid has a $k$-sparse generating set of its cycle space over $k$.
\end{obs}

\begin{proof}[Proof sketch.]
It suffices to prove this observation for connected matroids.
If $M$ is graphic, its atomic bonds form a sparse generating set. Conversely, if the dual matroid has a $k$-sparse generating set, build a graph whose vertices are the vectors in the sparse generating set.
\end{proof}

Following a similar route as in the plane, we obtain the following 3D MacLane Theorem.

\begin{cor}\label{mc}
Let $k$ be a field. A simply connected $k$-local 2-dimensional simplicial complex $C$ embeds in 3-space if and only if its cycle space $Z_2^k(C)$ has a $k$-sparse generating set.
\end{cor}

\begin{proof}
 Combine \autoref{MC-lem} and \autoref{intro-main}.
\end{proof}

\bibliographystyle{plain}
\bibliography{literatur}

\end{document}

%% file: minimalJ.tex
\usepackage{epsfig}
\usepackage{graphicx}
\usepackage{amsbsy}
\usepackage{amsmath}
\usepackage{amsfonts}
\usepackage{amssymb}
\usepackage{textcomp}
\usepackage{hyperref}
\usepackage{aliascnt}

\newcommand{\mcm}[3]{\newcommand{#1}[#2]{{\ensuremath{#3}}}} 

\mcm{\tuple}{1}{\langle #1 \rangle}
\mcm{\name}{1}{\ulcorner #1 \urcorner}
\mcm{\Nbb}{0}{\mathbb{N}}
\mcm{\Zbb}{0}{\mathbb{Z}}
\mcm{\Rbb}{0}{\mathbb{R}}
\mcm{\Cbb}{0}{\mathbb{C}}
\mcm{\Qbb}{0}{\mathbb{Q}}
\mcm{\Acal}{0}{\cal A}
\mcm{\Bcal}{0}{\cal B}
\mcm{\Ccal}{0}{\cal C}
\mcm{\Dcal}{0}{\cal D}
\mcm{\Ecal}{0}{\cal E}
\mcm{\Fcal}{0}{\cal F}
\mcm{\Gcal}{0}{\cal G}
\mcm{\Hcal}{0}{\cal H}
\mcm{\Ical}{0}{\cal I}
\mcm{\Jcal}{0}{\cal J}
\mcm{\Kcal}{0}{\cal K}
\mcm{\Lcal}{0}{\cal L}
\mcm{\Mcal}{0}{\cal M}
\mcm{\Ncal}{0}{\cal N}
\mcm{\Ocal}{0}{{\cal O}}
\mcm{\Pcal}{0}{{\cal P}}
\mcm{\Qcal}{0}{{\cal Q}}
\mcm{\Rcal}{0}{{\cal R}}
\mcm{\Scal}{0}{{\cal S}}
\mcm{\Tcal}{0}{{\cal T}}
\mcm{\Ucal}{0}{{\cal U}}
\mcm{\Vcal}{0}{{\cal V}}
\mcm{\Wcal}{0}{{\cal W}}
\mcm{\Xcal}{0}{{\cal X}}
\mcm{\Ycal}{0}{{\cal Y}}
\mcm{\Zcal}{0}{{\cal Z}}
\mcm{\Mfrak}{0}{\mathfrak M}

\mcm{\restric}{0}{\upharpoonright}
\mcm{\upset}{0}{\uparrow}
\mcm{\onto}{0}{\twoheadrightarrow}
\mcm{\smallNbb}{0}{{\small \mathbb{N}}}
\DeclareMathOperator{\preop}{op}
\mcm{\op}{0}{^{\preop}}

{\begin{array}{c}
\setlength{\unitlength}{1em}}%
{\end{array}}

\usepackage{amsthm}

\newcommand{\theoremize}[2]{\newaliascnt{#1}{thm} \newtheorem{#1}[#1]{#2} \aliascntresetthe{#1}}

\theoremstyle{plain}
\newtheorem{thm}{Theorem}[section]
\theoremize{lem}{Lemma}
\theoremize{skolem}{Skolem}
\theoremize{fact}{Fact}
\theoremize{sublem}{Sublemma}
\theoremize{claim}{Claim}
\theoremize{obs}{Observation}
\theoremize{prop}{Proposition}
\theoremize{cor}{Corollary}
\theoremize{que}{Question}
\theoremize{oque}{Open Question}
\theoremize{oprob}{Open Problem}
\theoremize{con}{Conjecture}

\theoremstyle{definition}
\theoremize{dfn}{Definition}
\theoremize{rem}{Remark}
\theoremize{eg}{Example}
\theoremize{exercise}{Exercise}
\theoremstyle{plain}